\documentclass[a4paper,12pt]{article}
\usepackage{comment}
\usepackage{cite}
\usepackage{amsmath}
\usepackage{amssymb}
\usepackage{amsfonts}
\usepackage[T1]{fontenc}
\usepackage[utf8]{inputenc}
\usepackage{graphicx}
\usepackage{fancyhdr}
\usepackage{float}
\usepackage{xcolor}
\usepackage{authblk}
\usepackage{mathrsfs}
\usepackage{empheq}
\usepackage[hyphens]{url}
\usepackage{hyperref} 
\usepackage[]{breakurl}

\usepackage{geometry}
\begin{document}

\title{From the Schaar and L\"osch-Schoblick integrals to representations of the Glaisher-Kinkelin constant}

\author[$\dagger$]{Jean-Christophe {\sc Pain}$^{1,2,}$\footnote{jean-christophe.pain@cea.fr}\\
\small
$^1$CEA, DAM, DIF, F-91297 Arpajon, France\\
$^2$Universit\'e Paris-Saclay, CEA, Laboratoire Mati\`ere en Conditions Extr\^emes,\\ 
F-91680 Bruy\`eres-le-Ch\^atel, France
}

\maketitle

\begin{abstract}
In this article, we present two integral representations of the logarithm of the Glaisher-Kinkelin constant, relying on two different integral formulations of the so-called Binet function $\mu(x)$. The first one is attributed to Schaar (and also often referred to as ``the second Binet formula''), and the second one is due to L\"osch and Schoblik. It seems that the two new expressions (formulas (28) and (33) of the present article) of the Glaisher-Kinkelin constant, can not be easily deduced from know ones.
\end{abstract}

\section{Introduction}

The Glaisher-Kinkelin constant $A$ \cite{Kinkelin1860,Glaisher1878,Finch2003} plays a major role in different fields of mathematics and theoretical physics. For instance, it appears in the Gaudin model of quantum spin chains \cite{Gaudin1976}. $A$ is also related to prime numbers, through \cite{Gorder2012}
\begin{equation}
\prod _{p{\text{ prime}}}p^{\frac {1}{p^{2}-1}}={\frac {A^{12}}{2\pi e^{\gamma }}},
\end{equation}
where $\gamma$ is the usual Euler-Mascheroni constant. The Glaisher-Kinkelin constant can be expressed as
\begin{equation}
    A=\lim _{n\rightarrow \infty }{\frac {\left(2\pi \right)^{\frac {n}{2}}n^{{\frac {n^{2}}{2}}-{\frac {1}{12}}}~e^{-{\frac {3n^{2}}{4}}+{\frac {1}{12}}}}{G(n+1)}},
\end{equation}
involving in the denominator the Barnes $G$-function: 
\begin{equation}
    G(n)=\prod_{k=1}^{n-2}k!.
\end{equation}
One has for instance
\begin{equation}
G(1/2)={\frac {2^{1/24}e^{1/8}}{A^{3/2}\pi ^{1/4}}}
\end{equation}
or also
\begin{equation}
G(1/4)={\frac {1}{2^{9/16}A^{9/8}\pi ^{3/16}\varpi ^{3/8}}}\exp \left({\frac {3}{32}}-{\frac {\beta(2)}{4\pi }}\right),
\end{equation}
where $\beta$ is the Dirichlet beta function\footnote{$\beta(2)$ is also the Catalan constant.}:
\begin{equation}
    \beta (s)=\sum _{n=0}^{\infty }\frac {(-1)^{n}}{(2n+1)^{s}},
\end{equation}
and 
\begin{equation}
\varpi = \frac {\Gamma (1/4)^{2}}{2{\sqrt {2\pi }}} 
\end{equation}
denotes the transcendental lemniscate constant (ratio of the perimeter of Bernoulli's lemniscate to its diameter).

The logarithm of the Glaisher-Kinkelin constant $A$ can be expressed in different manners, such as \cite{Glaisher1878}:
\begin{equation}\label{gla}
    \int_0^{1/2}\log\left[\Gamma(x+1)\right]~\mathrm{d}x=-\frac{1}{2}-\frac{7}{24}\log 2+\frac{1}{4}\log\pi+\frac{3}{2}\log A, 	
\end{equation}
where $\Gamma$ represents the usual Gamma function, or also \cite{Glaisher1878,Almkvist1998}:
\begin{equation}
    \log A=\frac{1}{12}-2\int_0^{\infty}\frac{x\log x}{e^{2\pi x}-1}~\mathrm{d}x.
\end{equation}
Other (but few) integral representations are actually available in the literature (see for instance the work of Choi and Nash \cite{Choi1997}. More recently, we obtained \cite{Pain2024a}:
\begin{align}
    \log A&=\frac{1}{9}\log 2+\frac{1}{24}+\frac{1}{3}\int_0^{\infty}\frac{\left(1-e^{-t/2}\right)\left[t~\coth\left(t/2\right)-2\right]}{t^3}~\mathrm{d}t,
\end{align}
\begin{equation}
    \log A=\frac{1}{3}+\frac{7}{36}\log 2-\frac{1}{6}\log\pi+\frac{1}{12}\int_0^{\infty}\frac{e^{-t}\left[(8-3t)~e^t-8~e^{t/2}-t\right]}{t^2(e^{t}-1)}~\mathrm{d}t, 	
\end{equation}
as well as \cite{Pain2024b}:
\begin{align}
    \log A=&\frac{1}{3}+\frac{7}{36}\log 2-\frac{\log\pi}{6}\nonumber\\
    &+\frac{2}{3}\int_0^{\infty}\left[\frac{e^{-t}}{8}-\frac{1}{(1+t)^{3/2}\log^2(1+t)}-\frac{1}{2}\frac{(\log(1+t)-2)}{(1+t)\log^2(1+t)}\right]~\frac{\mathrm{d}t}{t},  	
\end{align}
and
\begin{equation}
    \log A=\frac{\log 2}{36}+\frac{1}{3}\int_0^{\infty}\left[\frac{\tanh\left(\frac{t}{4}\right)}{t}-\frac{e^{-t}}{4}-\right]\,\mathrm{d}t.
\end{equation}

Using two different improper integral representations of the so-called Binet function \cite{Campbell1966}:
\begin{equation}\label{bin}
    \mu(x)=\log\left[\Gamma(x+1)\right]-\left(x+\frac{1}{2}\right)\,\log x+x-\log(\sqrt{2\pi}),
\end{equation}
the first one due to Schaar (see Sec. \ref{sec1}), and the second one to L\"osch and Schoblik (see Sec. \ref{sec2})), we derive two integral representations of the logarithm of the Glaisher-Kinkelin constant $\log A$, which, to our knowledge, were not published elsewhere.

\section{Integral representation deduced from the Schaar formula}\label{sec1} 

The following important result is due to Schaar \cite{Whittaker1990,Campbell1966,Gronwall1918} (a proof is given in Appendix A): 
\begin{equation}\label{mu1}
    \mu(x)=2x\int_0^{\infty}\frac{\arctan(t)}{e^{2\pi xt}-1}\,\mathrm{d}t.
\end{equation}
Making the change of variable $u=2\pi xt$, one gets
\begin{equation}
    \mu(x)=\frac{1}{\pi}\int_0^{\infty}\frac{\arctan\left(\frac{u}{2\pi x}\right)}{e^{u}-1}\,\mathrm{d}u.
\end{equation}
We have, integrating Eq. (\ref{bin}) from 0 to $1/2$:
\begin{equation}
    \int_0^{1/2}\mu(x)\,\mathrm{d}x= \int_0^{1/2}\log\left[\Gamma(x+1)\right]\,\mathrm{d}x+\int_0^{1/2}\left\{\left(x+\frac{1}{2}\right)\,\log x-x+\log(\sqrt{2\pi})\right\}\,\mathrm{d}x.
\end{equation}
It is easy to prove that
\begin{equation}
\int_0^{1/2}\left\{\left(x+\frac{1}{2}\right)\,\log x-x+\log(\sqrt{2\pi})\right\}\,\mathrm{d}x=\frac{1}{16}\left(-7-2\log 2+4\log\pi\right).
\end{equation}
Using the Fubini-Tonelli theorem, one can write
\begin{equation}
\int_0^{1/2}\mu(x)\,\mathrm{d}x=\frac{1}{\pi}\int_0^{1/2}\int_0^{\infty}\frac{\arctan\left(\frac{u}{2\pi x}\right)}{e^{u}-1}\,\mathrm{d}u\,\mathrm{d}x=\frac{1}{\pi}\int_0^{\infty}\int_0^{1/2}\frac{\arctan\left(\frac{u}{2\pi x}\right)}{e^{u}-1}\,\mathrm{d}x\,\mathrm{d}u
\end{equation}
or
\begin{equation}
\int_0^{1/2}\mu(x)\,\mathrm{d}x=\frac{1}{\pi}\int_0^{\infty}\frac{1}{e^{u}-1}\int_0^{1/2}\arctan\left(\frac{u}{2\pi x}\right)\,\mathrm{d}x\,\mathrm{d}u.
\end{equation}
One has to calculate the integral
\begin{equation}
\int_0^{1/2}\arctan\left(\frac{u}{2\pi x}\right)\,\mathrm{d}x.
\end{equation}
We have
\begin{equation}
  \arctan\left(x\right)+\arctan\left(\frac{1}{x}\right)=\frac{\pi}{2}  
\end{equation}
and thus
\begin{equation}
    \int_0^{1/2}\arctan\left(\frac{u}{2\pi x}\right)\,\mathrm{d}x=\int_0^{1/2}\left[\frac{\pi}{2}-\arctan\left(\frac{2\pi x}{u}\right)\right]\,\mathrm{d}x=\frac{\pi}{4}-\frac{u}{2\pi}\int_0^{\frac{\pi}{u}}\arctan(v)\,\mathrm{d}v.
\end{equation}
Integrating the last integral by parts gives
\begin{equation}
    \int_0^{1/2}\arctan\left(\frac{u}{2\pi x}\right)\,\mathrm{d}x=\frac{\pi}{4}-\frac{u}{2\pi}(v\,\arctan(v)\Big|_0^{\frac{\pi}{u}}+\frac{u}{2\pi}\int_0^{\frac{\pi}{u}}\frac{v}{v^2+1}\,\mathrm{d}v
\end{equation}
and thus
\begin{equation}
    \int_0^{1/2}\arctan\left(\frac{u}{2\pi x}\right)\,\mathrm{d}x=\frac{\pi}{4}-\frac{1}{2}\arctan\left(\frac{\pi}{u}\right)+\frac{u}{4\pi}\int_0^{\frac{\pi}{u}}\log(1+v^2)\,\mathrm{d}v.
\end{equation}
The latter result enables us to write
\begin{equation}
    \int_0^{1/2}\arctan\left(\frac{u}{2\pi x}\right)\,\mathrm{d}x=\frac{\pi}{4}-\frac{1}{2}\arctan\left(\frac{\pi}{u}\right)+\frac{u}{4\pi}\,\log\left(1+\frac{\pi^2}{u^2}\right).
\end{equation}
We thus have
\begin{align}
    \int_0^{1/2}\log\left[\Gamma(x+1)\right]\,\mathrm{d}x=&\frac{1}{\pi}\int_0^{\infty}\frac{1}{e^{u}-1}\left[\frac{\pi}{4}-\frac{1}{2}\arctan\left(\frac{\pi}{u}\right)+\frac{u}{4\pi}\,\log\left(1+\frac{\pi^2}{u^2}\right)\right]\,\mathrm{d}u\nonumber\\
    &+\frac{1}{16}\left(-7-2\log 2+4\log\pi\right).
\end{align}
Inserting this expression in Eq. (\ref{gla}) yields
\begin{align}\label{res1}
    \log A&=\frac{\log 2}{9}+\frac{1}{24}+\frac{2}{3\pi}\int_0^{\infty}\frac{1}{e^{t}-1}\left[\frac{\pi}{4}-\frac{1}{2}\arctan\left(\frac{\pi}{t}\right)+\frac{t}{4\pi}\,\log\left(1+\frac{\pi^2}{t^2}\right)\right]~\mathrm{d}t,
\end{align}
which is the first main result of the present work.

\section{Integral representation resulting from a representation due to L\"osch and Schoblik}\label{sec2}

L\"osch and Schoblik proposed the following integral representation of the Binet function \cite{Losch1951,Erdelyi1981}:
\begin{equation}\label{mu2}
    \mu(x)=-\frac{x}{\pi}\int_0^{\infty}\frac{\log\left(1-e^{-2\pi t}\right)}{t^2+x^2}\,\mathrm{d}t.
\end{equation}
Using the Fubini theorem, one can write
\begin{align}
\int_0^{1/2}\mu(x)\,\mathrm{d}x=&-\frac{1}{\pi}\int_0^{1/2}\int_0^{\infty}\frac{x\,\log\left(1-e^{-2\pi t}\right)}{t^2+x^2}\,\mathrm{d}t\,\mathrm{d}x\nonumber\\
=&-\frac{1}{2\pi}\int_0^{\infty}\log\left(1-e^{-2\pi t}\right)\int_0^{1/2}\frac{2x}{t^2+x^2}\,\mathrm{d}x\,\mathrm{d}t,
\end{align}
where
\begin{equation}
\int_0^{1/2}\frac{2x}{t^2+x^2}\,\mathrm{d}x=\log\left(1+\frac{1}{4t^2}\right).
\end{equation}
We thus have
\begin{align}
    \int_0^{1/2}\log\left[\Gamma(x+1)\right]~\mathrm{d}x&=-\frac{1}{2\pi}\int_0^{\infty}\log\left(1-e^{-2\pi t}\right)\log\left(1+\frac{1}{4t^2}\right)\,\mathrm{d}t\nonumber\\
    &+\frac{1}{16}(-7-2\log 2+4\log\pi),
\end{align}
which, inserted into Eq. (\ref{gla}) again, yields
\begin{equation}\label{res2}
    \log A=\frac{1}{24}+\frac{\log 2}{9}-\frac{1}{3\pi}\int_0^{\infty}\log\left(1-e^{-2\pi t}\right)\log\left(1+\frac{1}{4t^2}\right)\,\mathrm{d}t, 	
\end{equation}
which is the second main result of the present work. 

\section{Conclusion}

In this article, we proposed two integral representations of the logarithm of the Glaisher-Kinkelin constant. Both are based on a definite integral representation involving $\log\left[\Gamma(x+1)\right]$, which can itself be expressed as an improper integral, either using a representation originally proposed by Schaar and Binet or utilizing a formula by L\"osch and Schoblik. These two new expressions may prove useful in exploring new properties of the Glaisher-Kinkelin constant. Other relations should be obtained fro integral representations of the derivative of the Binet function (see Appendix B).

\appendix

\section{Derivation of the Schaar-Binet formula}

The Stieltjes constants are the coefficients of the Laurent expansion of the Riemann Zeta function $\zeta(x)$ around $x=1$ \cite{Connon2009}:
\begin{equation}
    \zeta(x)=\frac{1}{x-1}+\sum_{n=0}^{\infty}\frac{(-1)^n}{n!}\gamma_n(x-1)^n.
\end{equation}
In the case of the Hurwitz Zeta function $\zeta(x,z)$ \cite{Apostol1976}, we have the generalized Stieltjes constants:
\begin{equation}
    \zeta(x,z)=\sum_{n=0}^{\infty}\frac{1}{(n+z)^x}=\frac{1}{x-1}+\sum_{n=0}^{\infty}\frac{(-1)^n}{n!}\gamma_n(z)(x-1)^n.
\end{equation}
Let us differentiate $\zeta(x,z)$ with respect to variable $x$ and set $x=0$. Using the Lerch identity
\begin{equation}
    \log\Gamma(z)=\zeta'(0,z)+\frac{1}{2}\log(2\pi),
\end{equation}
one gets
\begin{equation}\label{star2}
    \frac{1}{2}\log(2\pi)-\log\Gamma(z)=1+\sum_{n=0}^{\infty}\frac{\gamma_{n+1}(z)}{n!}.
\end{equation}
Now, using the Coffey expression of the generalized Stieltjes constants
\cite{Coffey2007}:
\begin{equation}\label{cof}
    \gamma_n(z)=\frac{1}{2z}\log^nz-\frac{1}{n+1}\log^{n+1}z-2\Re\int_0^{\infty}\frac{i(z+iv)\log^n(z-iv)}{(z^2+v^2)(e^{2\pi v}-1)}\,\mathrm{d}v,
\end{equation}
and inserting expression (\ref{cof}) into Eq. (\ref{star2}) yields
\begin{align}
    \sum_{n=0}^{\infty}\frac{\gamma_{n+1}(z)}{n!}=&\frac{1}{2z}\sum_{n=0}^{\infty}\frac{1}{n!}\log^{n+1}z-\sum_{n=0}^{\infty}\frac{1}{n!}\frac{1}{(n+2)}\log^{n+2}z\nonumber\\
    &-2\Re\int_0^{\infty}\frac{i(z+iv)\sum_{n=0}^{\infty}\frac{1}{n!}\log^{n+1}(z-iv)}{(z^2+v^2)(e^{2\pi v}-1)}\,\mathrm{d}v.
\end{align}
Then, using a Taylor expansion series, one gets
\begin{equation}
    \sum_{n=0}^{\infty}\frac{1}{n!}\log^{n+1}z=z~\log z.
\end{equation}
One has also
\begin{equation}
    \sum_{n=0}^{\infty}\frac{y^{n+2}}{(n+2)!}=e^y-1-y,
\end{equation}
or equivalently
\begin{equation}
    \sum_{n=0}^{\infty}\frac{1}{n!}\frac{y^{n}}{(n+2)}=\frac{e^y(y-1)-1}{y^2},
\end{equation}
and thus
\begin{equation}
    \sum_{n=0}^{\infty}\frac{1}{n!}\frac{\log^{n+2}z}{(n+2)}=z(\log z-1)+1.
\end{equation}
Since
\begin{equation}
    \sum_{n=0}^{\infty}\frac{1}{n!}\log^{n+1}(z-iv)=(z-iv)\log(z-iv),
\end{equation}
one gets
\begin{equation}
    \sum_{n=0}^{\infty}\frac{\gamma_{n+1}(z)}{n!}=\frac{\log z}{2}-z(\log z-1)-1-2\int_0^{\infty}\frac{\arctan(\frac{x}{z})}{e^{2\pi x}-1}~\mathrm{d}x,
\end{equation}
yielding
\begin{equation}
    \mu(x)=2x\int_0^{\infty}\frac{\arctan(t)}{e^{2\pi xt}-1}\,\mathrm{d}t.
\end{equation}
This expression was also obtained by Ramanujan \cite{Ramanujan} for $x>0$.

\section{Expressions for the function $\mu(x)$ and its derivative}

In addition to expressions (\ref{mu1}) and (\ref{mu2}), other formulas are known for the Binet function $\mu(x)$, such as the following one, attributed to Binet as well \cite{Sasvari1999}:
\begin{equation}
    \mu(x)=\int_0^{\infty}e^{-xt}\left(\frac{1}{e^t-1}-\frac{1}{t}+\frac{1}{2}\right)~\frac{\mathrm{d}t}{t}.
\end{equation}
Some formulas are available also for the function $\nu(z)$ defined as
\begin{equation}
\nu(x)=\frac{d\mu(x)}{dx}+\frac{1}{2x}.
\end{equation} 
The two following ones are due to Poisson:
\begin{equation}
    \nu(x)=\int_0^{\infty}e^{-xt}\left(\frac{1}{e^t-1}-\frac{1}{t}\right)~\mathrm{d}t,
\end{equation}
and
\begin{equation}
    \nu(x)=\frac{1}{2x}+2\int_0^{\infty}\frac{t}{(1+t^2)(1-e^{2\pi xt})}~\mathrm{d}t.
\end{equation}

\end{document}